\documentclass[a4paper]{article}

\usepackage[utf8]{inputenc}
\usepackage{amsmath}
\usepackage{amssymb}
\usepackage{amsfonts}
\usepackage{mathtools}
\usepackage{graphicx}
\usepackage{tabularx}
\usepackage{tabularcalc}
\usepackage{booktabs}
\usepackage{multirow}
\usepackage[ruled]{algorithm2e}
\usepackage{subcaption}
\usepackage{siunitx}
\usepackage[breaklinks, colorlinks, linkcolor=blue, citecolor=blue, urlcolor=black]{hyperref}
\usepackage{doi}

\usepackage[numbers]{natbib}

\usepackage{xspace}
\usepackage{xcolor}

\usepackage{tikz,pgfplots}
\pgfplotsset{compat=1.16}

\newcommand{\MIP}{{MIP}\xspace}
\newcommand{\MIPs}{{MIPs}\xspace}
\newcommand{\MINLP}{{MINLP}\xspace}

\newcommand{\st}{\text{s.t.\xspace}}

\newcommand{\foralltext}{\text{for all\ }}

\newcommand{\testset}[1]{\textsc{#1}}

\title{A diving heuristic for mixed-integer problems with unbounded semi-continuous variables}

\newcommand{\myorcidlink}[1]{\,\href{https://orcid.org/#1}{\raisebox{-0.45ex}{\includegraphics[width=1.8ex]{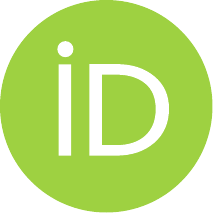}}}}
\makeatletter
\renewcommand*{\@fnsymbol}[1]{\ensuremath{\ifcase#1 \or * \or a\or b\or
		c \or d \or e \or f \or \dagger\dagger
		\or \ddagger\ddagger \else\@ctrerr\fi}}
\makeatother

\author{
	Katrin Halbig\protect\myorcidlink{0000-0002-8730-3447}\footnote{Corresponding author}
	\thanks{\parbox [t] {\linewidth} {Friedrich-Alexander-Universität Erlangen-Nürnberg, Department of Data Science,\newline Cauerstr.~11, 91058 Erlangen, Germany, katrin.halbig@fau.de}}
	\and Alexander Hoen\protect\myorcidlink{0000-0003-1065-1651}\thanks{Zuse Institute Berlin, Takustr.~7, 14195 Berlin, Germany, hoen@zib.de}
	\and Ambros Gleixner\protect\myorcidlink{0000-0003-0391-5903}$^b$\thanks{\parbox [t] {\linewidth} {Hochschule für Technik und Wirtschaft Berlin, 10313 Berlin, Germany,\newline gleixner@htw-berlin.de}}
	\and Jakob Witzig\protect\myorcidlink{0000-0003-2698-0767}\thanks{SAP SE, Dietmar-Hopp-Allee 17, 69190 Walldorf, Germany, jakob.witzig@sap.com}
	\and Dieter Weninger\protect\myorcidlink{0000-0002-1333-8591}\thanks{\parbox [t] {\linewidth} {Friedrich-Alexander-Universität Erlangen-Nürnberg, Department of Mathematics,\newline Cauerstr.~11, 91058 Erlangen, Germany, dieter.weninger@fau.de}}
}

\date{October 2024}

\begin{document}

\maketitle

\begin{abstract}
Semi-continuous decision variables arise naturally in many real-world applications.
They are defined to take either value zero or any value within a specified range, and occur mainly to prevent small nonzero values in the solution.
One particular challenge that can come with semi-continuous variables in practical models is that their upper bound may be large or even infinite.
In this article, we briefly discuss these challenges, and present a new diving heuristic tailored for mixed-integer optimization problems with general semi-continuous variables.
The heuristic is designed to work independently of whether the semi-continuous variables are bounded from above, and thus circumvents the specific difficulties that come with unbounded semi-continuous variables.
We conduct extensive computational experiments on three different test sets, integrating the heuristic in an open-source MIP solver.
The results indicate that this heuristic is a successful tool for finding high-quality solutions in negligible time.
At the root node the primal gap is reduced by an average of 5\,\% up to 21\,\%, and considering the overall performance improvement, the primal integral is reduced by 2\,\% to 17\,\% on average.
\end{abstract}



\let\symbols\undefined

\newcommand\mysymbol[3]{%
\protected\gdef#1{#2}%
\ifdefined\symbols
\item[$#2$]#3
\fi
}

\ifdefined\symbols
\section*{List of symbols}

\begin{description}
\fi
\mysymbol{\varother}{x}{other variables}
\mysymbol{\varsemi}{y}{semi-cont. variables}
\mysymbol{\varbin}{z}{binary indicator variables}
\mysymbol{\nvarsemi}{\nu}{number of semi-cont./indicator variables}
\mysymbol{\nvarother}{\eta}{number of other variables}
\mysymbol{\consmat}{A}{constraint matrix}
\mysymbol{\conscoef}{a}{constraint coefficient}
\mysymbol{\lhs}{b}{right-hand side vector}
\mysymbol{\lb}{\ell}{lower bound vector}
\mysymbol{\ub}{u}{upper bound vector}
\mysymbol{\objother}{c}{objective coefficient vector of \varother}
\mysymbol{\objsemi}{d}{objective coefficient vector of \varsemi}
\mysymbol{\objbin}{e}{objective coefficient vector of \varbin}
\mysymbol{\setintvar}{\mathcal{I}}{index set of integer variables}
\mysymbol{\setvarsemi}{\mathcal{J}}{index set of semi-continuous/integer variables}
\mysymbol{\idxvarsemi}{j}{index of semi-cont./indicator variables}
\mysymbol{\minlot}{s}{minimal lotsize of semi-continuous/integer variables}
\mysymbol{\round}{\sigma}{function for rounding direction}
\mysymbol{\bound}{\beta}{function for new value of bound}
\mysymbol{\score}{\psi}{score function}
\mysymbol{\setcand}{\mathcal{K}}{index set of candidate variables}
\mysymbol{\cand}{\kappa}{index of candidate variable}
\ifdefined\symbols
\item[$\ast^\top$] use \verb|\top| for transpose
\item[$\textforall$] use macro
\end{description}
\fi

\section{Introduction}
\label{sect:introduction}

In mathematical optimization, semi-continuous variables
arise naturally in models for many real-world applications, and have been studied since at least 1979~\cite{Land_Powell_1979}.
Formally, a variable $\varsemi$ is called \emph{semi-continuous} if it is defined to take either the value 0 or any value within the range specified by its \emph{lower semi-continuous bound} $\minlot$ and \emph{upper semi-continuous bound} $\ub$:
\begin{equation}
	\varsemi \in \{0\} \cup [\minlot,\ub] \text{ with } 0 < \minlot \leq \ub \text{ and } \ub \in \mathbb{R}_+ \cup \{\infty\}.
\end{equation}
This characteristic makes semi-continuous variables a powerful modeling tool for capturing scenarios where a variable's influence is either fully absent or continuous within a specific range.

Semi-continuous variables were introduced to simplify the solution of blending problems when materials must be excluded from the blend if they cannot
be used in significant quantities~\cite{Beale_1985}.
Nowadays, they are used in a wide range of real-world applications.
In portfolio optimization, semi-continuous variables are often used to avoid a large number of small trades or holdings which are undesirable because of management and transaction costs \cite{Perold1984}.
Also in unit commitment problems in electrical power systems, semi-continuous variables are used to model minimum online and offline times of units \cite{Chen2019}.
Supply chain management problems are also problems in which semi-continuous variables often occur. Here, semi-continuous variables are used to model economic order quantities, technical requirements 
because transportation of low quantities is undesirable from an operational point of view,
or to describe the state of a machine that is either turned off and thus produces nothing or turned on and thus has to produce at least a minimal amount (so-called \emph{minimum lot size}) due to economic and technical reasons; see \cite{Erlenkotter1990Ford, Harris1913Many, Park2015, Timpe2000}.

In the present article, we distinguish between two cases of semi-continuous variables.
If $\ub$ is finite, we say the variable is \emph{bounded semi-continuous} and if $\ub$ is infinite, we refer to it as an \emph{unbounded semi-continuous} variable.
In particular, we are interested in the common scenario where semi-continuous variables are part of a mixed-integer linear optimization problem.
Formally, we consider an optimization problem of the form
\begin{subequations}\label{eq:semicontMIP}
	\begin{align}
	\underset{\varother,\varsemi}{\min} \quad & \objother^\top \varother + \sum_{\idxvarsemi=1}^{\nvarsemi}f_\idxvarsemi(\varsemi_\idxvarsemi)\label{eq:semicontMIP:obj}\\
	\st \quad& \consmat (\varother,\varsemi) \geq \lhs,\label{eq:semicontMIP:consmat}\\
	& \varsemi_\idxvarsemi \in \{0\} \cup [\minlot_\idxvarsemi,\ub_\idxvarsemi], && \foralltext \idxvarsemi \in \setvarsemi,\label{eq:semicontMIP:semibound}\\
	& \varother_i \in \mathbb{Z}, && \foralltext i \in \setintvar_\varother,\label{eq:semicontMIP:integer_x}\\
	& \varsemi_i \in \mathbb{Z}, && \foralltext i \in \setintvar_\varsemi,\label{eq:semicontMIP:integer_y}
	\end{align}
\end{subequations}
where $y$ denotes a vector of semi-continuous variables with lower semi-continuous bound $\minlot \in \mathbb{R}_+^\nvarsemi$ and upper bound $\ub \in (\mathbb{R}_+ \cup \{\infty\})^\nvarsemi$.
Further variables are summarized in vector $x$.
The objective function is defined by vector
$\objother \in \mathbb{R}^\nvarother$, which relates to $\varother$, and by discontinuous functions
\begin{equation*}
        f_\idxvarsemi \colon \{0\} \cup [\minlot_\idxvarsemi,\ub_\idxvarsemi] \to \mathbb{R},\ 
	f_\idxvarsemi(\varsemi_\idxvarsemi)\coloneqq
	\begin{cases}
	0, &\text{ if } \varsemi_\idxvarsemi = 0,\\
	\objsemi_\idxvarsemi \varsemi_\idxvarsemi + \objbin_\idxvarsemi, &\text{ if } \varsemi_\idxvarsemi \geq \minlot_\idxvarsemi,
	\end{cases}
\end{equation*}
which relate to $\varsemi_\idxvarsemi$ for all $\idxvarsemi\in\setvarsemi \coloneqq \{1,\dots,\nvarsemi\}$.
Here, $\objsemi_\idxvarsemi \in \mathbb{R}$ models unit costs and
$\objbin_\idxvarsemi \in \mathbb{R}$ models set-up costs.
Constraint~\eqref{eq:semicontMIP:semibound}
specifies all $\nvarsemi$ semi-continuous variables.
Further linear constraints are given by a constraint matrix $\consmat \in \mathbb{R}^{m \times (\nvarother+\nvarsemi)}$ with right-hand side vector $\lhs \in \mathbb{R}^m$.
Sets $\setintvar_\varother \subseteq \{1,\dots,\nvarother\}$ and $\setintvar_\varsemi \subseteq \setvarsemi$ specify the integer variables.
Note that variables $\varsemi_\idxvarsemi$ with $\idxvarsemi \in \setintvar_\varsemi$ are also called \emph{semi-integer}, but for the sake of simplicity, we will not distinguish these cases.

In order to motivate that a heuristic approach is appropriate to solve problem~\eqref{eq:semicontMIP}, 
let us briefly discuss its complexity.
The optimization problem~\eqref{eq:semicontMIP} can be transformed into the related decision problem
\begin{equation}\label{eq:semicontMIP:decisionproblem}
\exists\ (\varother, \varsemi)
\text{ satisfying \eqref{eq:semicontMIP:consmat}-\eqref{eq:semicontMIP:integer_y} with value } \objother^\top \varother + \sum_{\idxvarsemi=1}^{\nvarsemi}f_\idxvarsemi(\varsemi_\idxvarsemi)  \leq k,\ k \in \mathbb{Z}.
\end{equation}
This decision problem~\eqref{eq:semicontMIP:decisionproblem} is
$\mathcal{NP}$-hard, because
the decision problem 0-1 integer program~(see \cite{Wolsey2020}) is polynomially reducible to~\eqref{eq:semicontMIP:decisionproblem} by not using any variables $\varother$ and requiring
$\minlot_\idxvarsemi = \ub_\idxvarsemi = 1$ and
$\objbin_\idxvarsemi =0$ for all
$\idxvarsemi \in \setvarsemi$.
Since the decision problem~\eqref{eq:semicontMIP:decisionproblem} is no more difficult to solve than the optimization problem~\eqref{eq:semicontMIP}, it is evident that problem~\eqref{eq:semicontMIP} is intractable even for finite $\ub$.

Solution approaches for mixed-integer problems with semi-continuous variables include, for example, strengthening of semi-continuous bounds~\cite{Achterberg_etal_2020,ZHANG2024107074}, branch-and-bound methods taking into account properties of semi-continuous variables~\cite{Beale_1985,Farias2004}, generation of semi-continuous cuts~\cite{Farias2004}, convex hull descriptions of problems with semi-continuous variables~\cite{Angulo2014,Farias_etal_2013}, and heuristics~\cite{Crevits_etal_2012,Hooshmand2023}. 
It is notable that, with the exception of~\cite{Angulo2014,Farias2004,Farias_etal_2013}, the vast majority of publications on semi-continuous variables is limited to the bounded case. 
An extensive literature survey of solution methods with a focus on nonlinear optimization problems with semi-continuous variables is provided in~\cite{Sun_etal_2013}.

Our contributions include, in Section~\ref{sect:lprelax}, a discussion of the particular challenges of solving MIPs with unbounded semi-continuous variables,
which arise when using an LP relaxation.
In Section~\ref{sect:algorithm} we propose a diving heuristic that is independent of the boundedness of the semi-continuous variables.
The peculiarity of this heuristic is that it performs a depth-first-search in the branch-and-bound tree that indirectly takes into account the characteristics of semi-continuous variables by treating indicator reformulations of the semi-continuous variables appropriately,
which means that the way the heuristic works can also be seen as a transfer of the branching approach from~\cite{Beale_1985,Farias2004}.
Moreover, a high-performance publicly available C implementation of this heuristic is provided in form of a tight integration with the open-source solver SCIP~\cite{BestuzhevaEtal2023}.
Computational results thereof are discussed in Section~\ref{sect:results}.
One of our test sets is comprised of real-world supply chain management instances, the other two test sets are publicly available and contain general \MIPs and artificially generated supply chain management instances.
A brief summary and a short discussion of open research questions for future work in Section~\ref{sect:conclusions} conclude the article.

\section{Formulations and LP Relaxation}
\label{sect:lprelax}

In this section, we will discuss formulations of semi-continuous variables in
problem~\eqref{eq:semicontMIP}, paying particular attention to the differences between bounded and unbounded semi-continuous variables.

Consider the feasible set $\{0\} \cup [\minlot,\ub]\subseteq \mathbb{R}$ of one 
semi-continuous variable.
This set is the union of two polyhedra.
For formulations of semi-continuous variables, it is now crucial how the convex hull of the union of the two polyhedra can be represented.
For $\ub = \infty$ the convex hull 
is not a \emph{bounded MIP-representable} set~\cite{Jeroslow1984,Vielma2015},
which means that no linear mixed-binary formulation exists.
If $\ub < \infty$, then the situation is less difficult since both polyhedra have the same recession cones and thus the convex hull of their union can itself be formulated as a polyhedron~\cite{Balas1988,Conforti_et_al_2010}.

If unbounded semi-continuous variables occur in problem~\eqref{eq:semicontMIP}, we are confronted with a more complex problem structure than if only bounded semi-continuous variables occur.
Therefore, the question arises as to whether it is necessary 
to consider a problem with unbounded semi-continuous variables at all.
In Section~\ref{sect:results} we substantiate that instances from practice actually contain unbounded semi-continuous variables.
In many cases, this is simply due to the fact that no feasible finite upper bound is known,
and determining a bound is of no interest to the practitioner since it would (a) not be a practical concern if $\varsemi$ takes a large value, and (b) it is known that unboundedness of $\varsemi$ does not render the entire problem unbounded.

Even more, it is sometimes not possible to determine a bound small enough to be used in numerical computations without knowing a near-optimal solution.
For example, in strategic supply chain planning, which involves optimizing operations over a period ranging from six months to multiple years, lot sizes are typically unbounded.
	Given the extended time horizon, time steps are discretized into daily, weekly, or monthly intervals.
	The feasibility of setting a maximum lot size (that is, upper bound $\ub$) largely depends on the nature of the goods being produced, transported, or procured, as well as their respective units of measurement.
	In some cases, expanding capacity to accommodate larger lot sizes can be straightforward, such as by leasing additional trucks.
	However, imposing finite maximum lot sizes can lead to numerical instability due to the substantial quantities involved within each time interval.
	Consequently, modelers often omit maximum lot sizes from their models when deemed unnecessary or impractical.
We therefore need to develop methods that can deal with unbounded semi-continuous variables \emph{during} the solving process.

There are several approaches to formulating semi-continuous variables.
Customary are the \emph{complementary formulation} resulting
in non-convex non-linear problems~\cite{Belotti2016,Stein_etal_2004} and \emph{disjunctive programming}
approaches~\cite{BALAS19793, Bonami2015}. 
The \emph{big-$M$-formulation}~\cite{Bonami2015, williams2013model} is presumably the easiest and commonly used disjunctive
programming approach, and is applicable 
in the case of a finite upper bound $\ub$.
By using a constant $M$ with $M \geq \ub$ a bounded
semi-continuous variable can be formulated with an additional binary variable $\varbin$ as
\begin{subequations}\label{eq:bigM}
	\begin{align}
	\varsemi \leq M \cdot \varbin,\label{eq:bigM:indcons} \\
	\minlot \cdot \varbin \leq \varsemi \leq \ub,\label{eq:bigM:lb} \\
	\varbin \in \{0,1\}.
	\end{align}
\end{subequations}
For $z = 1$ the constraint~\eqref{eq:bigM:indcons} is redundant and
$\minlot \leq  \varsemi \leq \ub$ is active.
In the other case $z = 0$, the variable $\varsemi$ is fixed to $0$ by~\eqref{eq:bigM:indcons} and~\eqref{eq:bigM:lb} whereby
$\varsemi \leq \ub$ is redundant.
In the case of an infinite upper bound $\ub$ one can use a
sufficiently large positive constant $M$, such that no optimal solution is cut off.

However, finding a suitable $M$ can be a major difficulty, see \cite{Kleinert2020} for a related discussion.
In practice, a pragmatic strategy is sometimes adopted by setting $M$ to a ``huge'' value, even though this may cut off optimal solutions.

Furthermore, the big-$M$~method may cause numerical problems.
When using a large value for $M$, 
the reciprocal $1/M$ can be close to the machine accuracy or the tolerances used by mathematical programming solvers when working with floating point arithmetic and become indistinguishable from zero.
This carries the risk of numerical difficulties during the solution process.
For example, if $\varbin$ takes a value close to zero, $M\cdot\varbin$ can still be large enough to set $\varsemi$ to a substantial value, whereas $\varbin$ is evaluated as integer feasible and zero.
In addition, the solutions of continuous relaxations may be weak, i.e., very far away from the optimal solution value.
Some of these difficulties can be defused by strengthening the $M$ value within the branch-and-bound tree by separation of local cuts~\cite{Chvatal_etal_2013} and bound propagation techniques~\cite{Achterberg_etal_2020,Savelsbergh_1994}.

Because of the above disadvantages
the focus of this article is on \emph{indicator constraints}~\cite{Belotti2016,Bonami2015},
which allows us to model bounded as well as unbounded semi-continuous variables exactly.
Indicator constraints activate a linear inequality based on the state of a 
binary variable known as the \emph{indicator variable}. 
When the indicator variable is set to a certain value, 
the corresponding constraint is included in the model; otherwise, it is ignored.
This representation is used henceforth and also in our later test sets, and has the advantage of being supported by virtually all state-of-the-art solvers.

To express a semi-continuous variable $\varsemi$ using indicator constraints, 
an additional binary indicator variable $\varbin$ is introduced. 
This indicator variable controls the constraint $\varsemi \leq 0$~\eqref{eq:indicator:indcons}. 
If $\varbin = 0$, it activates the constraint and forces the value of $\varsemi$ to be $0$.
To ensure that $\varsemi \geq \minlot$ holds when $\varbin = 1$,
an additional constraint $\varsemi \geq \minlot \cdot \varbin$ is introduced \eqref{eq:indicator:varbound}.
Overall, we obtain the formulation
\begin{subequations}\label{eq:indicator}
	\begin{align}
	\varbin = 0 \ \implies \ \varsemi \leq 0, \label{eq:indicator:indcons}\\
	\minlot \cdot \varbin \leq  \varsemi \leq \ub, \label{eq:indicator:varbound}\\
	\varbin \in \{0,1\}. \label{eq:indicator:bin}
	\end{align}
\end{subequations}
By reformulating all semi-continuous variables $\varsemi_\idxvarsemi$ for all $\idxvarsemi \in \setvarsemi$ in problem~\eqref{eq:semicontMIP} 
with approach \eqref{eq:indicator}, we arrive at
\begin{subequations}\label{eq:semicontMIP:ind}
	\begin{align}
	\underset{\varother,\varsemi,\varbin}{\min} \quad & \objother^\top \varother + d^\top \varsemi + e^\top z\label{eq:semicontMIP:obj:ind}\\
	\st \quad& \consmat (\varother,\varsemi) \geq \lhs,\label{eq:semicontMIP:consmat:ind}\\
	& \varbin_\idxvarsemi = 0 \ \implies \ \varsemi_\idxvarsemi \leq 0, && \foralltext \idxvarsemi \in \setvarsemi,\label{eq:semicontMIP:indicator:ind}\\
	& \minlot_\idxvarsemi \cdot z_\idxvarsemi \leq \varsemi_\idxvarsemi \leq u_j, && \foralltext \idxvarsemi \in \setvarsemi,\label{eq:semicontMIP:lowerbound:ind}\\
	& \varbin_\idxvarsemi \in \{0,1\}, && \foralltext \idxvarsemi \in \setvarsemi, \label{eq:semicontMIP:binary:ind}\\
	& \varother_i \in \mathbb{Z}, && \foralltext i \in \setintvar_\varother,\label{eq:semicontMIP:integer_x:ind}\\
	& \varsemi_i \in \mathbb{Z}, && \foralltext i \in \setintvar_\varsemi,\label{eq:semicontMIP:integer_y:ind}
	\end{align}
\end{subequations}
where $\objother$, $d$, $e$, $\consmat$, $\lhs$, $\setvarsemi$,
$\setintvar_\varother$, and $\setintvar_\varsemi$ are defined as in Section~\ref{sect:introduction}.

Branch-and-bound type algorithms~\cite{Dakin1965, LandDoig1960} 
rely on solving a relaxation to obtain information for bounding and branching decisions.
An LP relaxation is derived from a MIP by ignoring the integrality conditions on the integer variables.
If the LP relaxation is used to solve MIPs with a branch-and-bound enumeration approach, this is referred to as the \emph{LP-based branch-and-bound} method~\cite{Achterberg2009, Linderoth_Savelsbergh_1999}.
For an LP relaxation to a MIP involving indicator constraints as problem~\eqref{eq:semicontMIP:ind}, 
the integrality conditions~\eqref{eq:semicontMIP:binary:ind}, \eqref{eq:semicontMIP:integer_x:ind},
\eqref{eq:semicontMIP:integer_y:ind}, and the indicator constraint~\eqref{eq:semicontMIP:indicator:ind} are ignored.
Therefore, the semi-continuous variable $\varsemi$ can take any value in $[0,\ub]$, 
which often results in weak LP relaxations~\cite{Belotti2016}.
In particular, if objective coefficient $\objbin$ is positive, $\varbin$ can be set to zero to reduce set-up costs, but $\varsemi$ can be arbitrarily large in $[0,\ub]$.
We note that this LP relaxation is even worse than the LP relaxation of the big-$M$-formulation.
Due to numerical issues and the possible cutoff of optimal solutions, we nevertheless choose the formulation with indicator constraints in this article.

\section{A Tailored Diving Heuristic}
\label{sect:algorithm}

The determination of feasible solutions is essential for branch-and-bound type algorithms~\cite{Achterberg2009,Dakin1965,LandDoig1960}  in order to cut off parts of the enumeration tree and thus accelerate the solution process.
In Section~\ref{sect:introduction} we have  shown that problem~\eqref{eq:semicontMIP:decisionproblem} is $\mathcal{NP}$-Hard and thus it is not surprising that heuristics play a central role in solving problem~\eqref{eq:semicontMIP:ind}.

One class of heuristics are the so-called \emph{diving heuristics} (or \emph{dive-and-fix heuristics}).
Starting with a solution of a relaxation at the current node in a branch-and-bound tree diving heuristics strengthen variable bounds and reoptimize the relaxation iteratively.
This simulates a depth-first-search to a leaf in the branch-and-bound tree.
Diving heuristics for mixed-integer linear programs (\MIP) are described, for example, in~\cite{Berthold2008,WitzigGleixner2021,Wolsey2020},
and for mixed-integer nonlinear programs (\MINLP) in~\cite{Bonami2012}.

We now propose the diving heuristic \emph{Indicator Diving} (ID), which aims to determine a feasible solution for problem~\eqref{eq:semicontMIP:ind}.
Indicator Diving can get called at any node in a branch-and-bound tree after solving the LP relaxation  of problem~\eqref{eq:semicontMIP:ind} with solution $(\widehat{\varother},\widehat{\varsemi},\widehat{\varbin})$.
If the LP solution is feasible for~\eqref{eq:semicontMIP:ind}, there is nothing to do, and the heuristic is not called.
Otherwise, one has to specify a set $\setcand$ of candidate variables that are responsible for the infeasibility of $(\widehat{\varother},\widehat{\varsemi},\widehat{\varbin})$.
These are at first all integer variables $\varother_i$, $i \in \setintvar_\varother$, $\varsemi_i$, $i \in \setintvar_\varsemi$, and $\varbin_\idxvarsemi$, $\idxvarsemi \in \setvarsemi$, with fractional value.
In our specific case of \MIPs with indicator constraints, $\setcand$ additionally contains all binary indicator variables $\varbin_\idxvarsemi$,
if their corresponding indicator constraint is violated as well as
$\varbin_\idxvarsemi$ is not fixed already and $\widehat{\varbin}_\idxvarsemi$ is integer.
As long as the indicator variables are unfixed, the indicator constraints are omitted in the LP relaxation and thus the LP solution may violate these constraints.

For every candidate
$\cand \in \setcand \subseteq \setintvar_\varother \cup \setintvar_\varsemi \cup \setvarsemi$
three different functions get called:
\begin{itemize}
	\setlength{\itemsep}{0pt}
	\item a \emph{score} function $\score \colon \setcand \to \mathbb{R}$,
	\item a \emph{rounding} function $\round \colon \setcand \to \{\mathtt{up},\, \mathtt{down}\}$, and
	\item a \emph{bound value} function $\bound \colon \setcand \to \mathbb{R}$.
\end{itemize}
These functions are described in detail later in this section.

In Algorithm~\ref{alg:semicontMIP} the basic operation of Indicator Diving for \MIPs with semi-continuous variables as formulated in~\eqref{eq:semicontMIP:ind} is presented.
The diving heuristic repeatedly selects a candidate $\cand\in\setcand$ with maximal score $\score(\cand)$, and sets new bounds for the corresponding variable according to the results of $\round(\cand)$ and $\bound(\cand)$.
Afterwards, the LP gets optimized again with the updated bounds.
This step is optional and can be skipped at some or all iterations to keep the algorithm efficient.
Before updating the LP, it is also possible to propagate bound changes, that is, one tries to utilize the bound change of the candidate variable to strengthen further variable bounds.
After updating the LP, one can, in addition, try to find an integral solution via rounding~\cite{Achterberg2012}.
\begin{algorithm}
	\caption{Diving Heuristic for \MIPs with semi-cont. variables}
	\label{alg:semicontMIP}
	\KwIn{Problem~\eqref{eq:semicontMIP:ind} with LP relaxation,
		LP feasible solution $(\widehat{\varother},\widehat{\varsemi},\widehat{\varbin})$,
		score function $\score$, rounding function $\round$, bound value function $\bound$.}
	\KwOut{Feasible solution of~\eqref{eq:semicontMIP:ind} if one has been found.}
	$\begin{aligned}
		\setcand \coloneqq\ &
		{\{\idxvarsemi\in\setvarsemi \,\vert\, \widehat{\varbin}_\idxvarsemi \in \mathbb{Z},\ \varbin_\idxvarsemi \text{ unfixed},\ \eqref{eq:semicontMIP:indicator:ind} \text{ violated} \}}\\
		&\cup {\{i \in \setintvar_\varother \,\vert\, \widehat{\varother}_i \not\in \mathbb{Z}\}}
		\cup {\{i \in \setintvar_\varsemi \,\vert\, \widehat{\varsemi}_i \not\in \mathbb{Z}\}}
		\cup {\{\idxvarsemi \in \setvarsemi \,\vert\, \widehat{\varbin}_\idxvarsemi \not\in \mathbb{Z}\}}
	\end{aligned}$\\
	\While{$\setcand \not= \emptyset$}
	{
		\For{$\cand \in \setcand$}
		{
			Calculate score $\score(\cand)$, rounding direction $\round(\cand)$ and new value $\bound(\cand)$ of bound.
		}
		Select candidate $\cand$ with maximal score.\\
		$\setcand$ $\gets$ $\setcand \setminus \{\cand\}$\\
		\eIf{$\round(\cand) = \mathtt{up}$}
		{
			Increase lower bound to $\bound(\cand)$.
		}
		{
			Decrease upper bound to $\bound(\cand)$.
		}
		{(Optional)} propagate bound change.\\
		{(Optional)} update and solve LP.\\
		\eIf{LP infeasible}
		{
			Backtrack or abort.
		}
		{
			Update $\setcand$ based on new LP solution.
		}
		{(Optional)} round LP solution.\\
	}
\end{algorithm}

If the LP is infeasible, one can backtrack by undoing the last bound change and selecting the opposite rounding direction, or abort the heuristic without finding a feasible solution.
If the LP solution is feasible for the original problem~\eqref{eq:semicontMIP:ind}, the candidate set $\setcand$ is empty and the heuristic terminates with this feasible solution.

In order to describe the heuristic Indicator Diving completely, we define now the three mentioned functions $\score$, $\round$, and $\bound$.
For the score function $\score$ we have to distinguish between indicator variables and integer variables.
If we consider an indicator variable $\varbin_\cand$, the score is given by
\begin{equation}\label{eq:score:indicator}
	\score(\cand) \coloneqq
	\begin{cases}
	-1, &\text{ if } \widehat{\varsemi}_\cand \in \{0\} \cup [\minlot_\cand, \ub_\cand],\\
	100 \cdot (\minlot_\cand - \widehat{\varsemi}_\cand) / \minlot_\cand, &\text{ if } \widehat{\varsemi}_\cand \in (0,\minlot_\cand).
	\end{cases}
\end{equation}
For all other candidate variables we use another already existing diving strategy.
For example, one can switch to Farkas Diving \cite{WitzigGleixner2021}, as we will do for the computational tests in Section~\ref{sect:results}.
Since indicator variables should always be preferred, the score for other candidate variables is scaled to a range less than~$-1$.
For example, one can apply a sigmoid function, which has an ``S''-shaped graph, to map the scores to the interval $[-300,-100]$ while preserving the order, see~\cite{Halbig_DecHeur_2024, Schewe_Schmidt_Weninger:2019}.
As the candidate variable with the highest score is taken next, the indicator variables with semi-continuous variables outside their domain are treated first, followed by the indicator variables with semi-continuous variables inside their domain and thus actually already fulfilled, and finally, all other integer candidate variables are treated.

Moreover, the functions $\round$ and $\bound$ are defined as follows:
If candidate $\cand$ corresponds to an indicator variable $\varbin_\cand$, the rounding function is given by
\begin{equation}
\round(\cand) \coloneqq
\begin{cases}
\mathtt{up}, &\text{ if } \widehat{\varsemi}_\cand \geq 0.5\; \minlot_\cand,\\
\mathtt{down}, &\text{ if } \widehat{\varsemi}_\cand < 0.5\; \minlot_\cand,
\end{cases}
\end{equation}
and the bound value function is given by
\begin{equation}
\bound(\cand) \coloneqq
\begin{cases}
1, &\text{ if } \widehat{\varsemi}_\cand \geq 0.5\; \minlot_\cand,\\
0, &\text{ if } \widehat{\varsemi}_\cand < 0.5\; \minlot_\cand.
\end{cases}
\end{equation}
In other words, the binary indicator variable is fixed to one---and thus semi-continuous variable $\varsemi_\cand$ is fixed to $[\minlot_\cand,\ub_\cand]$---if the value of the corresponding semi-continuous variable in the LP solution is at least 50\,\% of the lower semi-continuous bound.
Otherwise, it is fixed to zero---and thus $\varsemi_\cand$ is also fixed to zero.
At this point, Indicator Diving differs from well-known diving heuristics in the literature, since the bounding step for one variable depends on the LP solution value of another variable.
For all other candidate variables, we switch to the same already existing diving heuristic as used for the score function.
Typically, such a diving heuristic uses for a candidate variable $\varother_\cand$ (or $\varsemi_\cand$, $\varbin_\cand$)
bound value $\bound(\cand) \coloneqq \lceil \widehat{\varother}_\cand \rceil$ if $\round(\cand) = \mathtt{up}$
and $\bound(\cand) \coloneqq \lfloor \widehat{\varother}_\cand \rfloor$ if $\round(\cand) = \mathtt{down}$.

\section{Computational Study}
\label{sect:results}

In this section, we present a comprehensive computational study of the diving heuristic proposed in Section~\ref{sect:algorithm}.
After describing the experimental setup, results for two main experiments are shown, for the first only the root node is processed and the second analyzes the overall performance impact.

\paragraph*{Computational setup.}
For the experiments, we used a pre-release version
of SCIP~9.0~\cite{bolusani2024} (git hash 7cc9c068bc) with
SoPlex~\cite{bolusani2024} version 6.0.3 (git hash 555f5d54) as underlying LP solver, running single-threaded.
The code is publicly available at~\url{https://github.com/scipopt}.
Indicator Diving is added as SCIP heuristic plugin \texttt{heur\_indicatordiving.c} using the generic diving algorithm framework of SCIP.
The output stream has been modified slightly in some cases to retrieve the desired information.

Roughly summarized, SCIP begins the solving process of a problem with a presolving phase to simplify the problem, for example, fixing variables or deleting redundant constraints.
The root node is then created and selected as the first open node.
After bound propagation is performed, the LP relaxation of the node is solved, potentially through multiple rounds in which cutting planes are added.
If the node is not solved to integer optimality, branching is performed to create new child nodes.
SCIP continues by selecting an open node and repeating this process until the problem is solved or predefined solving limits are reached.

During this process, primal heuristics are typically invoked at two key points for each node: before solving the LP relaxation and after solving the LP relaxation.
Since Indicator Diving requires an LP solution, it is called after solving the LP relaxation, and for our experiments, it is invoked only at the root node.

All presented computational results were generated on a compute cluster using compute nodes
with Intel Xeon Gold 6326 processors with 2.9 GHz and 32 GB RAM;
see \cite{WoodyClusterWebsite} for more details.

\paragraph*{Performance metrics.}

For the purpose of evaluating the performance of Indicator Diving, we compare \emph{shifted geometric means} of \emph{primal gaps} and \emph{primal integrals}.

For a single instance, let $\mathcal{R}_{max}$ be the total running time and $\mathcal{R}_1,\dots,\mathcal{R}_\tau \in [0,\mathcal{R}_{max}]$ be the points in time when a new incumbent solution is found, $\mathcal{R}_0\coloneqq 0$, $\mathcal{R}_{\tau+1}\coloneqq\mathcal{R}_{max}$.
Let $(\tilde{\varother},\tilde{\varsemi},\tilde{\varbin})_{opt}$ be an optimal or best known solution and let
$(\tilde{\varother},\tilde{\varsemi},\tilde{\varbin})_{\mathcal{R}}$ be the incumbent solution at point $\mathcal{R} \in [0,\mathcal{R}_{max}]$.

The primal gap function $\text{PG}\colon [0,\mathcal{R}_{max}] \to [0,1]$ is defined as
\begin{equation}
\text{PG}(\mathcal{R}) \coloneqq
\begin{cases}
1, \text{\qquad if no incumbent until point $\mathcal{R}$},\\[5pt]
0,  \text{\qquad if } | (\objother, \objsemi, \objbin)^\top (\tilde{\varother},\tilde{\varsemi},\tilde{\varbin})_{opt} | = | (\objother, \objsemi, \objbin)^\top (\tilde{\varother},\tilde{\varsemi},\tilde{\varbin})_{\mathcal{R}} | = 0,\\[5pt]
1,  \text{\qquad if } (\objother, \objsemi, \objbin)^\top (\tilde{\varother},\tilde{\varsemi},\tilde{\varbin})_{opt} \cdot (\objother, \objsemi, \objbin)^\top (\tilde{\varother},\tilde{\varsemi},\tilde{\varbin})_{\mathcal{R}} < 0,\\[5pt]
\displaystyle \frac{|(\objother, \objsemi, \objbin)^\top (\tilde{\varother},\tilde{\varsemi},\tilde{\varbin})_{opt} - (\objother, \objsemi, \objbin)^\top (\tilde{\varother},\tilde{\varsemi},\tilde{\varbin})_{\mathcal{R}}|} {\max\{| (\objother, \objsemi, \objbin)^\top (\tilde{\varother},\tilde{\varsemi},\tilde{\varbin})_{opt}|, |(\objother, \objsemi, \objbin)^\top (\tilde{\varother},\tilde{\varsemi},\tilde{\varbin})_{\mathcal{R}}|\}},
 \text{\quad else,}
\end{cases}
\end{equation}
and building on this the primal integral is defined as
\begin{equation}
\text{PI} \coloneqq \sum_{t=1}^{\tau + 1} \text{PG}(\mathcal{R}_{t-1}) \cdot (\mathcal{R}_{t} - \mathcal{R}_{t-1}).
\end{equation}
The primal integral
provides an absolute measurement for the performance of primal heuristics by considering the evolution of the incumbent solution over time, rewarding algorithms that find good solutions early.
For a detailed description of the primal integral see~\cite{Berthold2013}.
For the best-known solution of an instance as the base value we can use publicly available solutions (see \cite{MiplibWebsite}) and/or results of previous runs.

For average values over multiple instances, we use the \emph{shifted geometric mean} (SGM) of these values over all instances.
For values $w_1,\dots,w_N \geq 0$ and shift $s\geq0$ we determine the SGM by
\begin{equation}\label{eq:sgm}
\text{sgm}(w_1,\dots,w_N,s) \coloneqq \left(\prod_{i=1}^{N} (w_i + s)\right)^\frac{1}{N} - s.
\end{equation}
Unlike the arithmetic mean, the shifted geometric mean has the advantage that extreme values have less influence on the average.
For further discussion on evaluating computational results with SGMs see, for example,~\cite{Achterberg07a}.
If not stated otherwise, we use for all average values the shifted geometric mean with shift $s=1$.

\paragraph*{Test sets.}

We consider three sets of instances.
The first test set contains 588~real-world supply chain management instances (\testset{SCM}) with semi-continuous variables modeled with indicator constraints
supplied by our industry partner SAP~SE~\cite{SAP}.
The most important components are inventory holding, capacity restrictions, procurement,
transport, production, and demand fulfillment.
In~\cite{Gamrath2019} a more detailed description of very similar instances can be found.
The second test set contains 192 artificially generated supply chain management instances based on real-world supply chain management models, again with semi-continuous variables modeled by indicator constraints.
They represent a fictive company procuring components, producing cellphones of different types, transporting them to distribution centers, and satisfying customer demands.
This test set is provided by SAP~SE~\cite{SAP} and is publicly available at~\cite{cellphone}.
We refer to it hereinafter as \testset{Cellphone}.
The third set contains all 42 instances of the MIP\-LIB 2017~\cite{MIPLIB2017} Collection with indicator constraints, publicly available at~\cite{MiplibWebsite}.
We have not filtered these further, so they do not necessarily conform to our specifications, but necessary ones are checked when Indicator Diving is called.

Performance variability can be caused by the compute nodes as well as various random factors within SCIP.
For example, the score function~\eqref{eq:score:indicator} for indicator variables with a semi-continuous variable within its domain (which is constant $-1$) is also equipped with a random factor in our implementation.
This mechanism obtains the highest score (among a series of constants $-1$) regardless of the input sequence of the problem to SCIP.
To account for the effect of performance variability, we use three different seeds (including the default seed of zero) and treat every instance-seed combination as one individual observation.
From now on, we refer to such a combination simply as an instance, which triples the number of instances.
So, we have test set \testset{SCM} with 1764 instances, \testset{Cellphone} with 576 instances, and \testset{Miplib} with 126 instances.

In Table~\ref{tbl:statistics} some basic statistics of all three test sets are summarized.
We consider only instances that reached the calling point of ID, that is, instances for which the root LP was solved and ID was indeed called.
Grouped by test set, different numbers of the instances are reported, once from the original problem and once from the problem that ID has received.
In row ``vars'' and ``conss'' we count the number of all variables and all constraints, in row ``ind'' we count only the number of indicator constraints.
The latter number is further divided into the number of indicator constraints belonging to unbounded semi-continuous variables (``unbd'') and to bounded semi-continuous variables with very large upper bound ($10^4 < \ub < \infty$, ``big'').
The line was drawn at $10^4$, since SCIP adds an additional linear inequality $\varsemi \leq \ub \varbin$ for smaller bounds; see also \eqref{eq:bigM:indcons} for the big-M formulation.
In column ``total''
we state the sum of these counts over all instances per test set.
The minimum, maximum, and shifted geometric mean of these figures is reported in columns ``min'', ``max'', and ``sgm'', respectively. 

\begin{table}[ht]
	\begin{center}
	{\footnotesize \setlength\tabcolsep{5pt}
		\begin{tabular}{l@{\hspace{5pt}}l|rrrr|rrrr}
			\toprule
			& & \multicolumn{4}{c|}{original problem} & \multicolumn{4}{c}{received problem} \\[5pt]
			&       &     total &   min &     max & sgm &     total &   min &    max & sgm  \\
			\midrule
			\multirow{5}{*}{\rotatebox[origin=c]{90}{\testset{SCM}}}
			\multirow{5}{*}{\rotatebox[origin=c]{90}{(\#1300)}}	&	vars &  14420569 &  1034 &  594067 &   6090.87 &   8345699 &   345 &  229155 &   4119.30 \\
			 &	conss &   8481884 &   287 &  515321 &   2329.60 &   4007191 &   166 &  225346 &   1493.93 \\
			&	ind &    362939 &     7 &   29599 &     78.94 &    209101 &     2 &   24083 &     67.74 \\
			&	unbd &    315083 &     0 &   29599 &     76.75 &       276 &     0 &      46 &      0.02 \\
			&	big &     47856 &     0 &   23754 &      0.05 &    104519 &     0 &   23844 &      5.45 \\
			\midrule
			\multirow{5}{*}{\rotatebox[origin=c]{90}{\testset{Cellphone}}}
			\multirow{5}{*}{\rotatebox[origin=c]{90}{(\#463)}} &	vars & 6778402 &  7065 &  24123 &  13498.18 &   75880 &    32 &    493 &    124.13 \\
			 &	conss & 4119674 &  4407 &  14467 &   8241.05 &   67253 &    27 &    439 &    109.75 \\
			&	ind &  144149 &   187 &    415 &    298.43 &   15614 &     4 &    107 &     24.55 \\
			&	unbd &  144149 &   187 &    415 &    298.43 & 	  0 &     0 &      0 &      0.00 \\
			&	big & 	  0 &     0 &      0 &      0.00 &   14162 &     4 &     96 &     22.25 \\
			\midrule
			\multirow{5}{*}{\rotatebox[origin=c]{90}{\testset{Miplib}}}
			\multirow{5}{*}{\rotatebox[origin=c]{90}{(\#99)}} & vars &3417906 &  1923 &  315484 &  20765.78 & 2081992 &   626 &  304610 &   8797.25 \\
			 & conss &4564677 &  3912 &  268835 &  33283.11 & 2435727 &  2153 &  165383 &  16217.48 \\
			& ind & 525174 &    52 &   21247 &   2856.09 &  414986 &    24 &   46395 &   1883.76 \\
			& unbd &	  0 &     0 &       0 &      0.00 &	   0 &     0 &       0 &      0.00 \\
			& big &	  0 &     0 &       0 &      0.00 &	   0 &     0 &       0 &      0.00 \\
			\bottomrule
	\end{tabular}}
	\caption{Statistics on dimensions of the instances.}
	\label{tbl:statistics}
	\end{center}
\end{table}

As an example, for test set \testset{SCM} 1300 instances reached the calling point.
These instances have between 7 and 29599 indicator constraints in the original problem
with an average of 78.94 indicator constraints.
In the received problem the average number of indicator constraints reduces to 67.74.
It is noteworthy that in test set \testset{Cellphone} in the original problem all indicator constraints belong to unbounded semi-continuous variables, but in the received problem all indicator constraints belong to bounded semi-continuous variables, due to presolving processes which determine finite bounds.
Although there are no or only a few unbounded indicator constraints remaining in the test sets, indicator constraints are nonetheless a concern because many have still a big upper bound, leading to the issues discussed in Section~\ref{sect:lprelax}.
Note also that in test set \testset{Miplib} all indicator constraints 
are already bounded by a small bound in the original problem.
Therefore, strictly speaking, they do not exhibit the property of containing big or infinite upper bounds, which motivated the design of Indicator Diving.

\paragraph*{Performance results.}
We start the evaluation with a root node experiment, for which we
stop SCIP after the end of the root node or a time limit of 5\,h is reached.
Within this solution process, Indicator Diving was called, if possible, after solving the LP relaxation, along with other default heuristics of SCIP.
Calling ID is, for example, not possible if all indicator constraints were fixed or deleted due to presolving processes.

\begin{table}[ht]
	\begin{center}
	\begin{tabular}{lrrrrrrr}
		\toprule
		         & \multicolumn{3}{c}{\#instances} & \multicolumn{3}{c}{PG root} & \multirow{3}{2.3em}{\ \ PG \\\ heur}\\
		         \cmidrule(r){2-4}\cmidrule(lr){5-7}
		         &        &       & best  & with & w/o &       & \\
				 & called & found & found & ID   & ID   & ratio & \\
		\midrule
		\testset{SCM} (\#1764) & 1300 & 1052 & 287 & 0.15 & 0.19 & 0.79 &  0.30 \\
		\testset{Cellphone} (\#576) & 463 & 365 & 55 & 0.23 & 0.25 & 0.92 & 0.46  \\
		\testset{Miplib} (\#126) & 99 & 16 & 7 & 0.76 & 0.80 & 0.95 & 0.50 \\
		\bottomrule
	\end{tabular}
	\end{center}
	\caption{Root node experiment.}
	\label{tbl:root_node}
\end{table}
For test set \testset{SCM}, Indicator Diving was called on 1300 out of 1764 instances and found a feasible solution in 1052 cases.
Since diving heuristics in SCIP use an objective cutoff, all these solutions improve the primal bound at the point in time they are found.
In 287 cases the solution found by ID was still the best solution at the end of the root node.

These figures are summarized for all three test sets in Table~\ref{tbl:root_node}.
Moreover, we take a look at the primal gap after solving the root node (column ``PG root'') for all instances on which ID was called.
We calculate the shifted geometric mean of the relative primal gaps and compare it to a run without ID.
For \testset{SCM} the primal gap is reduced by 21\% on average, from 0.19 to 0.15.
The average reduction of the primal gap at the point in time when the solution of ID is passed to SCIP is given in the last column ``PG heur''. The latter figures take only instances on which ID found a solution into account.

Taking a look at Table~\ref{tbl:root_node}, we can conclude that on \testset{SCM} and \testset{Cellphone} Indicator Diving has a high probability to find an improving solution.
In addition, over all test sets ID yields a large reduction of up to 70\,\% of the primal gap at the point in time when the solution is passed to SCIP.
This reduction of the primal gap persists until the end of the root node, and is large, especially on \testset{SCM}, reaching up to 21\,\% on average.

In addition, we aim to verify that indicator constraints modeling semi-continuous variables are, in fact, problematic and that they are frequently violated by the LP relaxation solution.
To this end, for the two test sets \testset{SCM} and \testset{Cellphone}, we investigated how many indicator variables were actually diving candidates in the first iteration of Algorithm~\ref{alg:semicontMIP}.
For test set \testset{SCM}, 13\,\% of the indicator variables were diving candidates whereby 60\,\% of the associated semi-continuous variables were in the prohibited range $(0,\minlot)$.
Considering only semi-continuous variables that have a big or infinite upper bound, 8\,\%
of the related indicator variables were diving candidates and of these 50\,\%
of the semi-continuous variables were in the range $(0,\minlot)$.
For \testset{Cellphone}, 14\,\% of the indicator variables were diving candidates whereby 99\,\% of the associated semi-continuous variables were in $(0,\minlot)$.
Restricted to semi-continuous variables that have a big or infinite upper bound, the numbers are essentially the same.

\begin{figure}[ht]
	\captionsetup[subfigure]{justification=centering}
	\captionsetup[subfigure]{format=hang}
	\definecolor{colorrootplot}{rgb}{0.12,0.46,0.70}
	\pgfplotsset{
		root plot/.style = {
			tick align=outside,
			tick pos=left,
			ticklabel style = {font=\footnotesize},
			xlabel={\small PG with ID},
			xmin=-0.05, xmax=1.05,
			ymin=-0.05, ymax=1.05,
			unit vector ratio={1 1},
			scale=0.6
		},
		root addplot/.style = {
			draw=colorrootplot,
			fill=colorrootplot,
			mark=x,
			only marks,
			mark size=3.5pt,
			thick
		}
	}
	\begin{subfigure}[b]{0.345\textwidth}
		\begin{tikzpicture}
		\begin{axis}[root plot,
		ylabel={\small PG w/o ID}]
		\addplot[root addplot]
		table[x=pGap_master,y=pGap_no]{graphics/primal_gap_root_scm_data.txt};
		\end{axis}
		\end{tikzpicture}
		\caption{\testset{SCM} (\#1300)\\ Wins: 432\\ Losses: 154}
		\label{fig:primal_gap_root:scm}
	\end{subfigure}
	\begin{subfigure}[b]{0.29\textwidth}
		\begin{tikzpicture}
		\begin{axis}[root plot,
		yticklabels=\empty]
		\addplot[root addplot]
		table[x=pGap_master,y=pGap_no]{graphics/primal_gap_root_cellphone_data.txt};
		\end{axis}
		\end{tikzpicture}
		\caption{\testset{Cellphone} (\#463)\\ Wins: 106\\ Losses: 20}
		\label{fig:primal_gap_root:cellphone}
	\end{subfigure}
	\begin{subfigure}[b]{0.31\textwidth}
		\begin{tikzpicture}
		\begin{axis}[root plot,
		yticklabels=\empty]
		\addplot[root addplot]
		table[x=pGap_master,y=pGap_no]{graphics/primal_gap_root_miplib_data.txt};
		\end{axis}
		\end{tikzpicture}
		\caption{\testset{Miplib} (\#99)\\ Wins: 17\\ Losses: 8}
		\label{fig:primal_gap_root:miplib}
	\end{subfigure}
	\caption{Primal gap at end of root node computation.}
	\label{fig:primal_gap_root}
\end{figure}
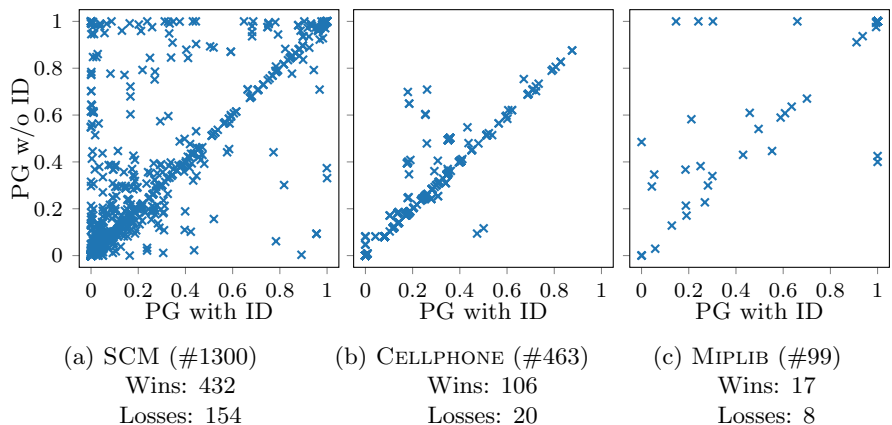

Furthermore, in Figure~\ref{fig:primal_gap_root} we compare the final primal gap achieved with and without ID over all instances for which ID was called.
The horizontal axis indicates the primal gap at the end of the root node computation for the run with ID, and the vertical axis indicates the same value in the run without ID.
Thus, crosses in the left upper corner represent instances with a reduction of the primal gap.
As crosses may overlap, we explicitly report the number of instances with a strict improvement (Wins) and a strict degradation (Losses) in the captions.
We consider an improvement/degradation to be strict if the absolute difference is at least $10^{-4}$.
These numbers demonstrate that there are considerably more instances with a strict improvement of the primal gap.
It is also noteworthy that within \testset{SCM} there are many instances where the gap decreases remarkably from a large value to almost zero due to Indicator Diving.

Finally, we analyze the overall performance impact of Indicator Diving on the complete tree search. For this we use a time limit of 20 minutes and no node limit. 
The run with Indicator Diving is compared to a run without Indicator Diving.
The results are summarized in Table~\ref{tbl:preformance} for each test set separately.
We disregard the instances where ID was not called at all, and consider the subsets of instances on which ID was called, and the subsets on which it was called and found a solution.
The numbers of instances of these subsets are stated in column ``\#inst'' and are slightly different compared to the numbers in the root node experiment above.
This is because we have different time limits and because ID can get called more than once due to restarts in SCIP.
We state the ratio of the shifted geometric means of the primal integrals of both runs in column ``ratio~PI'' whereby a value less than 1 indicates a performance improvement,
and the average time spent in ID (in seconds) is stated in column ``time~ID''.
Moreover, the average solving time of the run with ID and the average solving time of the run without ID is stated as well as the ratio thereof.
In the last two columns the number of instances that reached the time limit is specified.

\begin{table}[ht]
	\begin{center}
	{\small \setlength\tabcolsep{5pt}
		\begin{tabular}{llrrrrrrrrr}
			\toprule
			 & &       &     & & \multicolumn{3}{c}{solving time} &  \multicolumn{2}{c}{timeout}\\
			\cmidrule(r){6-8}\cmidrule(r){9-10}
			 & &        & ratio &  time   & with & w/o &       & with & w/o \\
			 & & \#inst & PI & ID &   ID &  ID & ratio &  ID  &  ID \\
			\midrule
			\multirow{2}{5.5em}{\testset{SCM} (\#1764)} & called & 1304  &  0.86  &  0.080 & 29.1 & 28.9 & 1.01  &  277  &  262\\
			& found & 1064  &  0.83  &  0.056 & 18.9 & 18.9 & 1.00  &  137  &  129 \\[5pt]
			\multirow{2}{5.5em}{\testset{Cellphone} (\#576)} & called & 463  &  0.97  &  0.004 & 161.1 & 177.0 &  0.91  &  281  &  289\\
			& found & 425  &  0.98  &  0.004 & 143.2 & 158.7 &  0.90 &  248  &  257 \\[5pt]
			\multirow{2}{5.5em}{\testset{Miplib} (\#126)} & called & 94  &  0.98  &  4.456 & 375.3 & 380.2 &  0.99 &  51  &  51 \\
			& found & 19  &  0.92  &  1.087 & 77.5 & 79.2 &  0.98  &  3  &  4\\
			\bottomrule
		\end{tabular}}
	\end{center}
	\caption{Performance comparison of SCIP with and without ID. Absolute times are given in seconds.}
	\label{tbl:preformance}
\end{table}

Taking a look at Table~\ref{tbl:preformance}, we can conclude that Indicator Diving yields a large reduction of the primal integral between 14\,\% and 17\,\% on test set \testset{SCM} and a reduction between 2\,\% and 8\,\% on the other two test sets.
Moreover, ID uses only a negligible proportion of the total solving time on all three test sets.
The impact on the overall solving time is almost neutral on \testset{SCM} and \testset{Miplib}, but reduced by up to 10\,\% on \testset{Cellphone}.
This aligns with previous findings that the impact of primal heuristics on the total solving time is minor~\cite{Berthold2013}.
Timeouts indicate that there is a considerable number of hard instances in our test sets even though the average solving time is moderate.
Their variation depends more on performance variability of the subsequent processes of SCIP than on ID.

To summarize, Indicator Diving helps the MIP solver to find better solutions earlier during the solving process, which is an important, if not the most important metric in practice.
As the total solving time is neutral or increases, we observe that its impact on performance is somewhat leveled out by the ensemble of other solving techniques applied by SCIP.

\section{Conclusion}
\label{sect:conclusions}

In this article, we discuss the challenges of using bounded and unbounded semi-continuous variables and propose a tailored diving heuristic 
for solving mixed-integer problems with semi-continuous variables that can be employed either standalone or integrated into a \MIP solver.
An implementation in C of this heuristic and two of the three treated
test sets are publicly accessible to be able to follow the extensive
computational experiments in detail.
One of the test sets was newly generated by our industry partner in the course of this research work and made available to the public in order to facilitate future research on the problems studied.

As part of this work, two computational experiments were carried out to evaluate the practical suitability of the diving heuristic.
The first one is a root node experiment, which exhibits a reduction of the primal gap on all three test sets when using the heuristic.
Particularly noteworthy is the result on one test set, where a reduction of the primal gap of 21\,\% is achieved.
In a second experiment we compared the performance of a MIP solver with and without the new heuristic
and observed that the use of the heuristic improves the primal integral by 2\,\% to 17\,\%.

Finally, we would like to address two research questions that complement the topic of the article presented.
First, we have found that unbounded semi-continuous variables occur in real-world instances and we have also shown in our computational study that unbounded semi-continuous variables are often transformed to bounded semi-continuous variables by bound propagation methods.
This leads to the assumption that specially adapted and more aggressively used domain reduction methods and presolving techniques might be able to transform even more unbounded semi-continuous variables into bounded semi-continuous variables.
Such approaches might also make it possible to determine tighter bounds for bounded semi-continuous variables.
Second, if an LP relaxation, in particular a simplex tableau, is used when generating cuts for problem~\eqref{eq:semicontMIP:ind}, then the indicator constraints are not taken into account.
Consequently, important information is lost in the cut generating process, which usually leads to weaker cuts. It would now be an interesting research question what possibilities there are to incorporate the information of the indicator constraints into the generation of valid cuts.

\section*{Acknowledgments}
\label{sec:acknowledgments}

The work for this article has been partly conducted within the Research Campus Modal funded by the German Federal Ministry of Education and Research (BMBF grant number 05M20ZBM).
The authors would like to thank SAP for its long-term support and for providing test instances.

\bibliographystyle{abbrvnat}
\bibliography{bibliography}

\end{document}